\documentstyle[12pt]{article}
\def\Ext{\mathop{\rm Ext}\nolimits}
\def\id{\mathop{\rm id}\nolimits}

\def\diag{\mathop{\rm diag}\nolimits}

\def\e{{\varepsilon}}

\def\ar{{\ \longrightarrow\ }}

\def\i{{\infty}}

\def\q{{\hfill $\Box$}}%\par\medskip}}
\def\proof{{\bf Proof. }}

\def\be{\begin{equation}}
\def\ee{\end{equation}}

\newtheorem{thm}{Theorem}[section]
\newtheorem{lem}[thm]{Lemma}

\newtheorem{cor}[thm]{Corollary}

\newtheorem{rmk}[thm]{Remark}
\newtheorem{prob}[thm]{Problem}

\newcommand{\norm}[1]{\| #1 \|}

\newcommand{\til}[1]{\widetilde{#1}}

\def\H{{\cal H}}
\def\K{{\cal K}}

\date{}
\author{V.~M.~Manuilov}
\title{Asymptotic homomorphisms into the Calkin algebra
\footnote{This research was partially supported by
RFBR (grant No 99-01-01201).}}

\begin{document}

\maketitle

 \begin{abstract}
Let $A$ be a separable $C^*$-algebra and let $B$ be a stable $C^*$-algebra
with a strictly positive element. We consider the (semi)group
$\Ext^{as}(A,B)$ (resp. $\Ext(A,B)$) of homotopy classes of asymptotic
(resp. of genuine) homomorphisms from $A$ to the corona algebra $M(B)/B$
and the natural map $i:\Ext(A,B)\ar\Ext^{as}(A,B)$. We show that if $A$ is
a suspension then $\Ext^{as}(A,B)$ coincides with $E$-theory of Connes and
Higson and the map $i$ is surjective. In particular any asymptotic
homomorphism from $SA$ to $M(B)/B$ is homotopic to some genuine homomorphism.

 \end{abstract}

%%%%%%%%%%%%%%%%%%%%%%%%%%%%%%%%%%%%%%%%%%%%%%%%%%%%%%%%%%%%%%%%%%%
{\protect\small\protect
\section{Introduction}

}

Let $A$, $B$ be $C^*$-algebras. Remind \cite{Connes-Higson} that a
collection of maps
 $$
\varphi=(\varphi_t)_{t\in[1,\i)}:A\ar B
 $$
is called an asymptotic homomorphism if for every $a\in A$ the map
$t\mapsto\varphi_t(a)$ is continuous and
if for any $a,b\in A$, $\lambda\in{\bf C}$ one has
 $$
\lim_{t\to\i}\norm{\varphi_t(ab)-\varphi_t(a)\varphi_t(b)}=0;
 $$
 $$
\lim_{t\to\i}\norm{\varphi_t(a+\lambda b)-\varphi_t(a)-\lambda\varphi_t(b)}=0;
 $$
 $$
\lim_{t\to\i}\norm{\varphi_t(a^*)-\varphi_t(a)^*}=0.
 $$
Two asymptotic homomorphisms $\varphi^{(0)}$ and $\varphi^{(1)}$ are
homotopic if there exists an asymptotic homomorphism $\Phi$ from $A$ to
$B\otimes C[0,1]$ such that its compositions with the evaluation maps at
$0$ and at $1$ coincide with $\varphi^{(0)}$ and $\varphi^{(1)}$
respectively. The set of homotopy classes of asymptotic homomorphisms
from $A$ to $B$ is denoted by $[[A,B]]$
\cite{Connes-Higson,Dadarlat-Loring}.

Throughout this paper we always assume that $A$
is {\it separable} and that $B$ has a strictly positive element and is
stable, $B\cong B\otimes\K$, where $\K$ denotes the $C^*$-algebra of
compact operators. We will sometimes write $B=B_1\otimes\K$, where
$B_1=B$, to distinguish $B$ from $B\otimes\K$ when necessary.

By $\Ext(A,B)$ we denote the set of homotopy classes of extensions of $A$ by
$B$. We identify extensions with homomorphisms into the Calkin
algebra $Q(B)=M(B)/B$ by the Busby invariant \cite{Busby}. Two extensions
$f_0,f_1:A\ar Q(B)$ are homotopic if there exists an extension
$F:A\ar Q(B\otimes C[0,1])$ such that its composition with the evaluation
maps at $0$ and at $1$ coincide with $f_0$ and $f_1$ respectively.

Similarly we denote by $\Ext^{as}(A,B)$ the set of homotopy classes of
asymptotic homomorphisms from $A$ to $Q(B)$. Two asymptotic homomorphisms
 $$
\varphi^{(i)}=(\varphi_t^{(i)})_{t\in[1,\i)}:A\ar Q(B),
\quad i=0,1,
 $$
are homotopic if there exists an asymptotic homomorphism
$\Phi=(\Phi_t)_{t\in[1,\i)}:A\ar Q(B\otimes C[0,1])$ such that its
compositions with the evaluation maps at $0$ and at $1$ coincide with
$\varphi^{(0)}$ and $\varphi^{(1)}$ respectively. Asymptotic homomorphisms
into $Q(B)$ are sometimes called {\it asymptotic extensions}.

All these sets are equipped with a natural group structure when $A$ is a
suspension, i.e. $A=SD=C_0({\bf R})\otimes D$ for some $C^*$-algebra $D$.

As every genuine homomorphism can be viewed as an asymptotic one,
so we have a natural map
 \be\label{i}
i:\Ext(A,B)\ar\Ext^{as}(A,B).
 \ee
It is well-known that usually there is much more asymptotic homomorphisms
than genuine ones, e.g. for $A=C_0({\bf R}^2)$ all
genuine homomorphisms of $A$ into $\K$ are homotopy trivial though
the group $[[C_0({\bf R}^2),\K]]$ coincides with $K_0(C_0({\bf
R}^2))={\bf Z}$ via the Bott isomorphism.

The main purpose of the paper is to prove epimorphity of the map (\ref{i})
when $A$ is a suspension. This makes a contrast with the case of
mappings into the compacts. As a by-product we get
another description of the $E$-theory in terms of asymptotic extensions.

The main tool in this paper is the Connes--Higson map \cite{Connes-Higson}
 $$
CH:\Ext(A,B)\ar [[SA,B]],
 $$
which plays an important role in $E$-theory.
Remind that for $f\in\Ext(A,B)$
this map is defined by
$CH(f)=(\varphi_t)_{t\in[1,\i)}$, where $\varphi$ is given by
 $$
\varphi_t:\alpha\otimes a\longmapsto \alpha(u_t)f'(a), \qquad
a\in A, \ \alpha\in C_0(0,1).
 $$
Here $f':A\ar M(B)$ is a set-theoretic lifting for
$f:A\ar Q(B)$ and $u_t\in B$ is a quasicentral approximate unit
\cite{Arveson} for $f'(A)$. We are going to show that by fine tuning of
this quasicentral approximate unit one can define also a map
 $$
\til{CH}:\Ext^{as}(A,B)\ar [[SA,B]]
 $$
extending $CH$ and completing the commutative triangle diagram
 $$
\begin{array}{ccc}
\Ext(A,B)\!\!\!\!\!&\stackrel{i}{\ar}&\!\!\!\!\!\Ext^{as}(A,B)\\
&\!\!\!\!\!\!\!\!{\scriptstyle CH}\searrow&
\downarrow{\scriptstyle \til{CH}}\\
&&\!\!\!\!\![[SA,B]].
\end{array}
 $$

We will show that the map $\til{CH}$ is a isomorphism when $A$ is a
suspension.

I am grateful to K.~Thomsen for his hospitality during my visit
to \AA rhus University in 1999 when the present paper was conceived.
%I am also grateful to A. S. Mishchenko for fruitful discussions.

%%%%%%%%%%%%%%%%%%%%%%%%%%%%%%%%%%%%%%%%%%%%%%%%%%%%%%%%%%%%%%%%%%%
{\protect\small\protect
\section{An extension of the Connes--Higson map}

}

A useful tool for working with asymptotic homomorphisms is the
possibility of discretization suggested in
\cite{Mish-Noor,Man-Mish,Man-Th2}.
Let $\Ext^{as}_{discr}(A,B)$ denote the set of homotopy classes of
discrete asymptotic homomorphisms $\varphi=(\varphi_n)_{n\in{\bf N}}:A\ar
Q(B)$ with the additional crucial property suggested by Mishchenko: for
every $a\in A$ one has
 \be\label{n+1}
\lim_{n\to\i}\norm{\varphi_{n+1}(a)-\varphi_n(a)}=0.
 \ee
In a similar way we define a set $[[A,B]]_{discr}$ as a set of homotopy
classes of discrete asymptotic homomorphisms with the property
(\ref{n+1}).

\begin{lem}\label{cont=discr}
One has
$[[A,B]]=[[A,B]]_{discr}$, \ $\Ext^{as}(A,B)=\Ext^{as}_{discr}(A,B)$.

\end{lem}

\proof
The first equality is proved in \cite{Man-Th1}. The second one can be
proved in the same way. For an asymptotic
homomorphism $\varphi=(\varphi_t)_{t\in[1,\i)}:A\ar Q(B)$ one can find an
infinite sequence of points $\{t_i\}_{i\in{\bf N}}\subset[1,\i)$
satisfying the following properties
\begin{enumerate}
\item
the sequence $\{t_i\}_{i\in{\bf N}}$ is non-decreasing and approaches
infinity;
\item
for every $a\in A$ one has \quad
 ${\displaystyle
\lim_{i\to\i}\sup_{t\in[t_i,t_{i+1}]}\norm{\varphi_t(a)-\varphi_{t_i}(a)}=0.
}
 $
\end{enumerate}
Then $\phi=(\varphi_{t_i})_{i\in{\bf N}}$ is a discrete asymptotic
homomorphism. It is easy to see that two homotopic asymptotic homomorphisms
define homotopic asymptotic homomorphisms and that two discretizations
$\{t_i\}_{i\in{\bf N}}$ and $\{t'_i\}_{i\in{\bf N}}$ satisfying the above
properties define homotopic discrete asymptotic homomorphisms too, hence
the map $\Ext^{as}(A,B)\ar\Ext^{as}_{discr}(A,B)$ is well defined. The
inverse map is given by linear interpolation of discrete asymptotic
homomorphisms.
\q

Let $(\varphi_n)_{n\in{\bf N}}$ be a discrete asymptotic homomorphism and
let $(m_n)_{n\in{\bf N}}$ be a sequence of numbers $m_n\in{\bf N}$. Then
we call the sequence
 $$
(\underbrace{\varphi_1,\ldots,\varphi_1}_{m_1\ {\rm times}},
\underbrace{\varphi_2,\ldots,\varphi_2}_{m_2\ {\rm times}},
\varphi_3,\ldots)
 $$
a {\it reparametrization} of the sequence $(\varphi_n)_{n\in{\bf N}}$.
It is easy to see that any reparametrization does not change the homotopy
class of an asymptotic homomorphism.

\begin{lem}\label{boundedness}
There exists a sequence of liftings $\varphi'_n:A\ar M(B)$ for
$\varphi_n$, which is continuous uniformly in $n$.

\end{lem}

\proof
It is easy to see \cite{Loring} that
 $$
\lim_{n\to\i}\sup_{n\leq k<\i}\norm{\varphi_k(a)-\varphi_k(b)}\leq
\norm{a-b}
 $$
for any $a,b\in A$.
By the Bartle--Graves selection theorem \cite{Bartle-Graves}, cf.
\cite{Loring} there exists a continuous selection $s:Q(B)\ar M(B)$. Put
$\varphi'_n(a)=s\varphi_n(a)$, $a\in A$.
\q

Now we are going to construct the map
$\til{CH}:\Ext^{as}(A,B)\ar[[SA,B]]$.
Due to Lemma \ref{cont=discr} it is
sufficient to define the map $\til{CH}$ as a map from
$\Ext^{as}_{discr}(A,B)$ to $[[SA,B]]_{discr}$.

For $a,b\in A$, $\lambda\in{\bf C}$ put
\begin{eqnarray*}
P_n(a,b)&=&\varphi_n(a)\varphi_n(b)-\varphi_n(ab);\\
L_n(a,b,\lambda)&=&\varphi_n(a)+\lambda\varphi_n(b)
-\varphi_n(a+\lambda b);\\
A_n(a)&=&\varphi_n(a)^*-\varphi_n(a^*)
\end{eqnarray*}
and define $P'_n(a,b)$, $L'_n(\lambda,a)$, $A'_n(a)$ in the
same way but with the liftings $\varphi'_n$ instead of $\varphi_n$.

In what follows we identify $B=B_1\otimes\K$ (resp. $M(B)$) with the
$C^*$-algebra of compact (resp. adjointable) operators on the standard
Hilbert $C^*$-module $B_1\otimes l_2({\bf N})=l_2(B_1)$ and use the notion
of diagonal operators in $B$ and $M(B)$ in this sense. The following Lemma
shows how one has to choose a quasicentral approximate unit that makes it
possible to define the map $\til{CH}$.

\begin{lem}\label{qcau}
Let $(\varphi_n)_{n\in{\bf N}}:A\ar Q(B_1\otimes\K)$ be a discrete
asymptotic homomorphism. Then there exists a reparametrization
of $(\varphi_n)_{n\in{\bf N}}$ and an approximate unit
$(u_n)_{n\in{\bf N}}\subset B_1\otimes\K$  with the following properties:
\begin{enumerate}
\item
for any $a\in A$ one has
 $$
\lim_{n\to\i}\norm{[\varphi'_n(a),u_n]}=0;
 $$
\item
for any $\alpha\in C_0(0,1)$, for any $a,b\in A$, $\lambda\in{\bf C}$
one has
 $$
\lim_{n\to\i}\norm{\alpha(u_n)P'_n(a,b)}=
\lim_{n\to\i}\norm{\alpha(u_n)L'_n(a,b,\lambda)}=
\lim_{n\to\i}\norm{\alpha(u_n)A'_n(a)}=0;
 $$
\item
$\lim_{n\to\i}\norm{u_{n+1}-u_n}=0$;
\item
every $u_n$ is a diagonal operator, $u_n=\diag\{u_n^1,u_n^2,\ldots\}$,
where diagonal entries $u_n^i$ belong to $B_1$ and
 $$
\lim_{i\to\i}\sup_n \norm{u_n^{i+1}-u_n^i}=0.
 $$
\end{enumerate}

\end{lem}

\proof
Let $\{F_n\}_{n\in{\bf N}}$ be a generating system for $A$
\cite{Connes-Higson}. This means that every $F_n\subset A$ is
compact,
\ $\ldots\subset F_n\subset F_{n+1}\subset\ldots$\ , \ $\cup_n F_n$ is
dense in $A$ and one has
 $$
F_n\cdot F_n\subset F_{n+m(n)};
\quad
F_n+\lambda F_n\subset F_{n+m(n)},  \ \ (|\lambda|\leq 1);
\quad
F_n^*\subset F_{n+m(n)}
 $$
for some integer sequence $m=(m_n)_{n\in{\bf N}}$.
Let also $\alpha_0=e^{2\pi ix}-1\in C_0(0,1)\cong C_0({\bf R})$ be a
(multiplicative) generator for $C_0({\bf R})$.

Put
 $$
\e_{n,k}=\sup_{a,b\in F_k,|\lambda|\leq 1}
\max(\norm{P_n(a,b)},\norm{L_n(a,b,\lambda)},\norm{A_n(a)}).
 $$
For every fixed $a,b,\lambda$ the sequences $(P_n(a,b))$,
$(L_n(a,b,\lambda))$ and $(A_n(a))$ vanish as $n$
approaches infinity, but the sequence $(\e_{n,n})_{n\in{\bf N}}$
does not have to vanish. Nevertheless
one can reparametrize the sequence $\{F_n\}$ by a sequence
$k=(k_n)_{n\in{\bf N}}$, which approaches infinity slowly enough
and such that $\e_{n,k(n)}$ vanishes as $n\to\i$. Put $\e_n=\e_{n,k(n)}$.
Then
 \be\label{eps}
\lim_{n\to\i}\e_n=0.
 \ee

Let $e=(e_n)_{n\in{\bf N}}\subset B$ be an approximate unit in $B$ and
let $Conv(e)$ denote its convex hull.

By induction we can
choose $u_n\in Conv(e)$ in such a way that $u_n\geq u_{n-1}$
and that the estimates
 \be\label{1}
\norm{[\varphi'_n(a),u_n]}<\e_n;
 \ee
and
 \be\label{2}
\norm{\alpha_0(u_n)P'_n(a,b)}<2\e_n;\quad
\norm{\alpha_0(u_n)L'_n(a,b,\lambda)}<2\e_n;\quad
\norm{\alpha_0(u_n)A'_n(a)}<2\e_n
 \ee
hold for any $a,b\in F_n$ and any $|\lambda|\leq 1$.

It is easy to see that the conditions (\ref{1}-\ref{2}) together with
Lemma \ref{boundedness} ensure the first two items of Lemma \ref{qcau}.

The above choice of $(u_n)_{n\in{\bf N}}$ does not yet ensure the condition
$\lim_{n\to\i}\norm{u_{n+1}-u_n}=0$. To make it hold we have
to renumber the sequence $(\varphi_n)_{n\in{\bf N}}$. At first divide every
segment $[u_n,u_{n+1}]$ into $n$ equal segments
$[u_{n_i},u_{n_{i+1}}]$, $i=1,\ldots,n$. Then as $0\leq u_i\leq 1$ for
all $i$, so we get $\norm{u_{n_{i+1}}-u_{n_i}}\leq\frac{1}{n}$. Finally we
have to change the sequences $(\varphi_1,\varphi_2,\varphi_3,\ldots)$ and
$(u_1,u_2,u_3,\ldots)$ by the sequence
$(\varphi_1,\varphi_2,\varphi_2,\varphi_3,\ldots)$, where each $\varphi_n$
is repeated $n$ times, and by the sequence
$(u_{1_1},u_{2_1},u_{2_2},u_{3_1},u_{3_2},u_{3_3},u_{4_1},\ldots)$
respectively.

To prove the last item of Lemma \ref{qcau} remind that
an approximate unit $(e_n)_{n\in{\bf N}}\in B=B_1\otimes\K$ can be
chosen to be diagonal, $e_n=b_n\otimes\epsilon_n$, where
$(b_n)_{n\in{\bf N}}\subset B_1$ and
$(\epsilon_n)_{n\in{\bf N}}\subset\K$ are approximate units in $B_1$
and in $\K$ respectively,
so the quasicentral approximate unit $(u_n)_{n\in{\bf N}}\subset
Conv(e)$ can be made diagonal as well, with diagonal entries from $B_1$.

Let $T$ be the right shift on the standard Hilbert $C^*$-module
$l_2(B_1)=B_1\otimes l_2({\bf N})$, $T\in M(\K)\subset M(B_1\otimes\K)$.
We can join $S$ to the
sets $\varphi'_n(F_n)$ in (\ref{1}) when constructing the sequence
$(u_n)$. Then the sequence $[T,u_n]\in B_1\otimes\K$
would vanish as $n$ approaches infinity. Hence
 \be\label{est11}
\lim_{n\to\i}\sup_i \norm{u_n^{i+1}-u_n^i}=0
 \ee
and the operators
 $$
\diag\{u_n^2-u_n^1,u_n^3-u_n^2,u_n^4-u_n^3,\ldots\}
 $$
are compact, so $\lim_{i\to\i}\norm{u_n^{i+1}-u_n^i}=0$.
Take $\e>0$. By (\ref{est11}) there exists some $N$ such that for all $n>N$
one has $\sup_i\norm{u_n^{i+1}-u_n^i}<\e$. Now consider the finite number of
compact operators $\diag\{u_n^2-u_n^1,u_n^3-u_n^2,u_n^4-u_n^3,\ldots\}$,
$1\leq n\leq N$. Due to their compactness there exists some $I$ such that
for $i>I$ one has $\norm{u_n^{i+1}-u_n^i}<\e$ for $1\leq n\leq N$.
Therefore for $i>I$ we have $\norm{u_n^{i+1}-u_n^i}<\e$ for every $n$,
i.e. $\sup_n\norm{u_n^{i+1}-u_n^i}<\e$.
\q

Put now
 $$
\til{CH}(\varphi)_n(\alpha\otimes a)=\alpha(u_n)\varphi'_n(a),
\qquad \alpha\in C_0(0,1), \ a\in A,
 $$
where $(u_n)_{n\in{\bf N}}$ satisfies the conditions of Lemma \ref{qcau}.
Items $i)-iii)$ of Lemma \ref{qcau} ensure that
$(\til{CH}(\varphi)_n)_{n\in{\bf N}}$
is a discrete asymptotic homomorphism from $SA$ to $B$. If
$(u_n)_{n\in{\bf N}}$ and $(v_n)_{n\in{\bf N}}$ are two quasicentral
approximate unities satisfying Lemma \ref{qcau} then the linear homotopy
$(tu_n+(1-t)v_n)_{n\in{\bf N}}$
provides that the maps $\til{CH}$ defined using these
approximate unities are homotopic. Finally, if $\varphi$ and $\psi$
represent the same homotopy class in $\Ext^{as}_{discr}(A,B)$ then
$\til{CH}(\varphi)$ and $\til{CH}(\psi)$ are homotopic.
If all $\varphi_n$ are constant, $\varphi_n=f:A\ar Q(B)$ with $f$ being a
genuine homomorphism, then obviously $CH(f)=\til{CH}(\varphi)$, so we have
\begin{lem}
The map $\til{CH}:\Ext^{as}(A,B)\ar [[SA,B]]$ is well defined and the
diagram
 $$
\begin{array}{ccc}
\Ext(A,B)\!\!\!\!\!&\stackrel{i}{\ar}&\!\!\!\!\!\Ext^{as}(A,B)\\
&\!\!\!\!\!\!\!\!{\scriptstyle CH}\searrow&
\downarrow{\scriptstyle \til{CH}}\\
&&\!\!\!\!\![[SA,B]].
\end{array}
 $$
is commutative. \q

\end{lem}

%%%%%%%%%%%%%%%%%%%%%%%%%%%%%%%%%%%%%%%%%%%%%%%%%%%%%%%%%%%%%%%%%%%
{\protect\small\protect
\section{An inverse for $\til{CH}$}

}

Let $\alpha_0=e^{2\pi ix}-1$ be a generator for $C_0(0,1)$ and let $T$ be
the right shift on the Hilbert space $l_2({\bf N})$. By
$q:M(B)\ar Q(B)$ we denote the quotient map. Define a
homomorphism
 $$
g:C_0(0,1)\ar Q(\K) \quad{\rm by}\quad
g(\alpha_0)=q(T)-1.
 $$

Remind that $B$ is stable and denote by $\iota:Q(B)\otimes\K\subset Q(B)$
the standard inclusion. Put
 $$
j=\iota\circ(g\otimes\id_B):SB\ar Q(\K)\otimes B\subset Q(B).
  $$
The homomorphism $j$ obviously induces a map
 $$
j_*:[[A,SB]]\ar \Ext^{as}(A,B).
 $$
Let $S:[[A,B]]\ar[[SA,SB]]$ denote the suspension map.
Then the composition $M=j_*\circ S$ gives a map
 $$
M:[[A,B]]\ar \Ext^{as}(SA,B).
 $$

Let
 \be\label{beta}
\beta=(\beta_n)_{n\in{\bf N}}:C_0({\bf R}^2)\ar \K
 \ee
be a discrete asymptotic homomorphism representing a generator of
$[[C_0({\bf R}^2),\K]]$.
For a discrete asymptotic extension $\varphi=(\varphi_n)_{n\in{\bf
N}}:A\ar Q(B)$ consider its tensor product by $\beta$
 $$
\varphi\otimes\beta=
(\varphi_n\otimes\beta_n)_{n\in{\bf N}}:S^2A\ar Q(B)\otimes\K
 $$
and denote its composition with the standard inclusion
$Q(B)\otimes\K\subset Q(B)$ by
 $$
Bott_1=\iota\circ(\varphi\otimes\beta):
\Ext^{as}(A,B)\ar\Ext^{as}(S^2A,B).
 $$
In a similar way define a map
 $$
Bott_2:[[A,B]]\ar [[S^2A,B]].
 $$

\begin{thm}\label{identity}
One has
 $$
M\circ\til{CH}=Bott_1;\quad
\til{CH}\circ M=Bott_2.
 $$
\end{thm}

\proof
We start with $M\circ\til{CH}=Bott_1$. Let $H$ be the standard Hilbert
$C^*$-module over $B$, $H=B\otimes l_2({\bf N})$. Put $\H=\oplus_{n\in
{\bf N}}H_n$, where every $H_n$ is a copy of $H$. We identify the
$C^*$-algebra of compact (resp. adjointable) operators on both $H$ and
$\H$ with $B$ (resp. $M(B)$). Instead of writing formulas in $Q(B)$ we
will write them in $M(B)$ and understand them modulo compacts.

Let $\varphi=(\varphi_n)_{n\in{\bf N}}:A\ar Q(B)$ represent an element
$[\varphi]\in\Ext^{as}_{discr}(A,B)$ and let $\varphi'_n:A\ar M(B)$ be
liftings for $\varphi_n$ as in Lemma \ref{boundedness}.

If $a_n:H_n\ar H_n$ is a sequence of operators then we write
$(a_1\oplus a_2\oplus a_3\oplus\ldots)$ for their direct sum acting on
$\H=\oplus_{n\in{\bf N}}H_n$. In what follows we use a shortcut
 $$
\alpha(u_n)\varphi'_n(a)=a_n.
 $$
Let $T$ be the right shift on $\H$, $T:H_n\ar H_{n+1}$.

Remind that $\alpha_0$ is a generator for $C_0(0,1)$ and that it is
sufficient to define asymptotic homomorphisms on the elements of the form
$\alpha\otimes a\otimes\alpha_0\in S^2A$.

The composition map $M\circ\til{CH}(\varphi)_n:S^2A\ar Q(B)$ acts by
 $$
M\circ\til{CH}(\varphi)_n(\alpha\otimes a\otimes\alpha_0)
=\Bigl(a_n\oplus a_n\oplus a_n\oplus\ldots\Bigr)(T-1)
 $$
modulo compacts on $\H$.

Let
 $$
v_n=\Bigl({\bf v}_n^1\oplus {\bf v}_n^2\oplus {\bf v}_n^3
\oplus\ldots\Bigr)\in M(B_1\otimes\K)
 $$
and
 $$
\lambda_n=\Bigl(\lambda_n^1\oplus \lambda_n^2\oplus \lambda_n^3
\oplus\ldots\Bigr)\in M(B_1\otimes\K)
 $$
be a direct sum of diagonal operators
${\bf v}_n^i=\diag\{v_n^i,v_n^i,v_n^i,\ldots\}$, $v_n^i\in B_1$,
and a direct sum of scalar operators,
$\lambda_n^i\in{\bf R}$,
(${\bf v}_n^i$ and $\lambda_n^i$ act on $H_i$).
Let the numbers $\lambda_n^i$ satisfy the properties
\begin{enumerate}
\item
$\lambda_n^1=0$ and $\lim_{i\to\i}\lambda_n^i=1$ for every $n$;
\item
$\lim_{n\to\i}\sup_i |\lambda_n^{i+1}-\lambda_n^i|=0$;
\item
$\lim_{n\to\i}\sup_i |\lambda_{n+1}^i-\lambda_n^i|=0$.
\end{enumerate}
We assume that the elements $v_n^i$ are selfadjoint and satisfy the
following properties:
\begin{enumerate}
\item
$\lim_{n\to\i}\sup_i\norm{v_n^{i+1}-v_n^i}=0$;
\item
$\lim_{n\to\i}\sup_i\norm{v_{n+1}^i-v_n^i}=0$;
\item
there exists a set of scalars $\lambda_n^i\in{\bf R}$, $n,i\in{\bf N}$,
satisfying the conditions $i)-iii)$ above and such that
 \begin{equation}\label{iii}
\lim_{n\to\i}\norm{(v_n^i-\lambda_n^i)b}=0
 \end{equation}
for every $i$ and for every $b\in B_1$.
\end{enumerate}

Let $p$ be a projection onto the first coordinate in
$H=B_1\otimes l_2({\bf N})$ and let $P=(p\oplus p\oplus p\oplus\ldots)$.
Then $P\lambda_n^i=\diag\{\lambda_n^i,0,0,\ldots\}$ and
the map $\beta_n$ (\ref{beta}) can be written as
 $$
\beta_n(\alpha\otimes\alpha_0)=P\cdot\alpha(\lambda_n)\cdot(T-1)
\in 1_{B_1}\otimes\K\subset M(B_1\otimes\K)
 $$
and the map $Bott_1(\varphi):S^2A\ar Q(B)$ can be written in the form
 $$
(Bott_1(\varphi))_n(\alpha\otimes a\otimes\alpha_0)=
\Bigl(\alpha(\lambda_n^1)\varphi'_n(a)\oplus
\alpha(\lambda_n^2)\varphi'_n(a)\oplus
\alpha(\lambda_n^3)\varphi'_n(a)\oplus \ldots\Bigr)
(T-1).
 $$
Consider also the path of asymptotic homomorphisms
$(\Phi_n(t))_{n\in{\bf N}}$, $t\in[0,1]$, given
by the formula
 $$
\Phi_n(t)(\alpha\otimes a\otimes\alpha_0)=
\Bigl(\alpha({\bf v}_n^1(t))\varphi'_n(a)\oplus
\alpha({\bf v}_n^2(t))\varphi'_n(a)\oplus
\alpha({\bf v}_n^3(t))\varphi'_n(a)\oplus \ldots\Bigr)
(T-1),
 $$
where for every $i$ $v_n^i(t)$ is a piecewise linear path connecting
$v_n^i(\frac{1}{k})=v_{n-1+k}^i$, $k\in{\bf N}$, and
$v_n^i(0)=\lambda_n^i$.
In view of (\ref{iii}) it is easy to check that $\Phi_n(t)$ is a homotopy
connecting the asymptotic homomorphisms $(Bott_1(\varphi))_n$ and
$\Phi_n=\Phi_n(0)$.

One of the obvious choices for $v_n^i$ is to put $(v_n^i)^{i\in{\bf
N}}=(0,\frac{1}{n},\frac{2}{n},\ldots,1,1,\ldots)$.
But for our purposes
it is better to use another choice. We take $v_n^i=u_i^n$ for all $n$ and
$i$.
Lemma \ref{qcau} ensures that the properties $i)-iii)$ are satisfied.

Now we have to connect the asymptotic homomorphisms
$(Bott_1(\varphi)_n)_{n\in{\bf N}}$ and
$(M\circ\til{CH}(\varphi)_n)_{n\in{\bf N}}$
by a homotopy in the class of asymptotic
homomorphisms. In fact we are going to do more and to connect
each of these asymptotic homomorphisms with a genuine homomorphism
$f:S^2A\ar Q(B)$ defined modulo compacts by
 $$
f(\alpha\otimes a\otimes\alpha_0)=\Bigl(
a_1\oplus a_2\oplus a_3\oplus\ldots
\Bigr)\cdot(T-1),
 $$
$\alpha\in C_0(0,1)$, $a\in A$.
Lemma \ref{qcau} ensures that $f$ is indeed a homomorphism.

At first we connect $f$ with $(M\circ\til{CH}(\varphi))_n$ by a path
$F_n(t)$, $t\in[0,1]$. Let $F_n(1)=(M\circ\til{CH}(\varphi))_n$.
Denote $\alpha(u_n)\varphi'_n(a)$ by $a_n$ and put
(modulo compacts)
 $$
F_n\left(\frac{1}{2}\right)
(\alpha\otimes a\otimes\alpha_0)=
\Bigl(
\underbrace{a_n\oplus
\ldots\oplus a_n}_{n\ {\rm times}}\oplus
a_{n+1}\oplus a_{n+1}\oplus a_{n+1}\ldots \Bigr)(T-1),
 $$
 $$
F_n\left(\frac{1}{3}\right)
(\alpha\otimes a\otimes\alpha_0)=
\Bigl(
\underbrace{a_n\oplus
\ldots\oplus a_n}_{n\ {\rm times}}\oplus
a_{n+1}\oplus a_{n+2}
\oplus a_{n+2}\oplus\ldots \Bigr)(T-1),
 $$
etc. Finally put
 $$
F_n(0)(\alpha\otimes a\otimes\alpha_0)=
\Bigl(
\underbrace{a_n\oplus
\ldots\oplus a_n}_{n\ {\rm times}}\oplus
a_{n+1}\oplus a_{n+2}
\oplus a_{n+3}\oplus\ldots \Bigr)(T-1)
 $$
and connect $F_n(1)$, $F_n(\frac{1}{2})$, $F_n(\frac{1}{3})$, $\ldots$ and
$F_n(0)$ by a piecewise linear path $F_n(t)$, $t\in[0,1]$.

It is easy to see
that for every $t>0$ the sequence $(F_n(t))_{n\in{\bf N}}$ is an
asymptotic homomorphism. And as the maps $F_n(0)$ and $f$ differ by
compacts, so they coincide as homomorphisms into $Q(B)$. Continuity in $t$
is also easy to check. So the asymptotic homomorphism
$(M\circ\til{CH}(\varphi)_n)_{n\in{\bf N}}$ is homotopic to the
homomorphism $f$.

Now we are going to construct a homotopy $F'_n(t)$, $t\in[0,1]$,
which connects $f$ with $(Bott_1(\varphi)_n)_{n\in{\bf N}}$.

For each $k\in{\bf N}$ consider the following sequence $({\bf
u}_n^k)_{n\in{\bf N}}$ of diagonal operators, each of which acts on the
corresponding copy of $H=H_n$ in their direct sum $\H$:
 $$
{\bf u}_n^k=\diag\{u_n^1,u_n^2,\ldots,u_n^{k-1},u_n^k,u_n^k,u_n^k\ldots\}.
 $$
Put $u_n(\frac{1}{k})={\bf u}_n^k$, $u_n(0)=u_n$ and connect them by a
piecewise linear path $u_n(t)$, $t\in[0,1]$. Then we get a strictly
continuous path of operators $u_n(t)$, which gives a homotopy
 $$
F'_n(t)(\alpha\otimes a\otimes\alpha_0)
=\Bigl(\underbrace{a_{1,n}(t)\oplus\ldots\oplus
a_{n,n}(t)}_{n\ {\rm times}}\oplus
a_{n+1,n+1}(t)\oplus
a_{n+2,n+2}(t)\oplus\ldots
\Bigr)(T-1),
 $$
where $a_{i,n}(t)=\alpha(u_i(t))\varphi'_n(a)$.

As
 $$
\Phi_n(\alpha\otimes a\otimes\alpha_0)
=\Bigl(a_{1,n}(1)\oplus a_{2,n}(1)\oplus a_{3,n}(1)\oplus\ldots
\Bigr)(T-1),
 $$
so for every $\alpha\otimes a$ one has
 $$
\lim_{n\to\i}\norm{F'_n(1)(\alpha\otimes a\otimes\alpha_0)-
\Phi_n(\alpha\otimes a\otimes\alpha_0)}=0,
 $$
hence the asymptotic homomorphisms $F'_n(1)$ and $\Phi_n$ are equivalent.
But we already know that $\Phi_n$ is homotopic to $Bott_1(\varphi)_n$.
On the other hand, it is easy to see that $F'_n(0)$
coincides with $f$ modulo compacts, so we can finally conclude that
$M\circ\til{CH}=Bott_1$ up to homotopy.

\smallskip
The second identity of Theorem \ref{identity} is much simpler to prove.
For $\psi=(\psi_n)_{n\in{\bf N}}:A\ar B$ we have (modulo compacts)
 $$
M(\psi)_n(\alpha_0\otimes a)=\Bigl(\psi_1(a)\oplus\psi_2(a)
\oplus\psi_3(a)\oplus\ldots\Bigr)(T-1),\qquad a\in A.
 $$
But as every $\psi_n(a)\in B_1\otimes\K$, i.e. is compact,
so when choosing a quasicentral
approximate unit $\{w_n\}_{n\in{\bf N}}$ for the map $(M(\psi))_{n\in{\bf
N}}$ we can define it by
 $$
w_n=\Bigl({\bf w}_1^{(n)}\oplus{\bf w}_2^{(n)}
\oplus{\bf w}_3^{(n)}\oplus\ldots\Bigr),
 $$
where each ${\bf w}_i^{(n)}$ is a finite rank diagonal operator of the form
 $$
{\bf w}_i^{(n)}=\diag\{\underbrace{\lambda_i b_n,
\ldots,\lambda_i b_n}_{m_n\ {\rm times}},0,0,\ldots\}
 $$
for some numbers $(m_n)_{n\in{\bf N}}$, where $(b_n)_{n\in{\bf N}}$ is
a quasicentral
approximate unit for $B_1$ and the scalars $\lambda_i$ are defined by
 $$
\lambda_i=
\left\lbrace\begin{array}{ll}
\frac{n-i+1}{n} & {\rm for}\  i<n,\\
\lambda_i=0 & {\rm for}\ i\geq n.
\end{array}\right.
 $$
But after such a choice of $w_n$ the map $(\til{CH}\circ M)(\psi)_n$
differs from the map $Bott_2(\psi)_n$ only by the presence of $b_n$, hence
these maps are equivalent.
\q

%%%%%%%%%%%%%%%%%%%%%%%%%%%%%%%%%%%%%%%%%%%%%%%%%%%%%%%%%%%%%%%%%%
{\protect\small\protect
\section{Case of $A$ being a suspension}

}

As there exists a homomorphism $C_0({\bf R})\ar C_0({\bf R}^3)\otimes
M_2$ that induces an isomorphism in $K$-theory and an asymptotic homomorphism
$C_0({\bf R}^3)\ar C_0({\bf R})\otimes\K$, which are inverse to each
other, so the groups $[[SA,B]]$ and $[[S^3A,B]]$  are
naturally isomorphic to each other and the same is true for the
groups $\Ext^{as}(SA,B)$ and $\Ext^{as}(S^3A,B)$. Hence we obtain
\begin{cor}\label{cor1}
If $A$ is a suspension then the map $\til{CH}:\Ext^{as}(A,B)\ar[[SA,B]]$
is an isomorphism. \q

\end{cor}

It was proved in \cite{Man-Th2} that if $A$ is a suspension
then the map
 $$
CH:\Ext(A,B)\ar[[SA,B]]
 $$
is surjective and the group $[[SA,B]]$ is contained in $\Ext(A,B)$ as a
direct summand. Hence from Corollary \ref{cor1} we immediately obtain
\begin{cor}\label{cor2}
Let $A$ be a suspension. Then
\begin{enumerate}
\item
the map $i:\Ext(A,B)\ar\Ext^{as}(A,B)$ is surjective,
hence every asymptotic extension $\varphi=(\varphi_t)_{t\in[1,\i)}:A\ar Q(B)$
is homotopic to a genuine extension;
\item
the group $\Ext^{as}(A,B)$ is contained in $\Ext(A,B)$ as a direct
summand. \q
\end{enumerate}

\end{cor}

\begin{prob}
{\rm
Is Corollary \ref{cor2} true when $A$ is not a suspension ?

}
\end{prob}

For $C^*$-algebras $A$ and $B$ consider the set of all extensions $f:A\ar
Q(B)$ that are homotopy trivial as asymptotic homomorphisms and denote by
$\Ext^{ph}(A,B)$ the set of homotopy classes of such homomorphisms. As
usual this set becomes a group when $A$ is a suspension.
We call the elements of $\Ext^{ph}(A,B)$ {\it phantom} extensions because
they constitute the part in $\Ext(A,B)$ which vanishes under the
suspension map $S:\Ext(A,B)\ar\Ext(SA,SB)$, cf.~\cite{HL-T}.

\begin{cor}
If $A$ is a second suspension then there is a natural decomposition
 $$
\Ext(A,B)=\Ext^{ph}(A,B)\oplus\Ext^{as}(A,B).
 $$

\vspace{-2.5em}
\q

\end{cor}

\begin{rmk}
{\rm
If $A$ is both a nuclear $C^*$-algebra and a suspension then the groups
$\Ext(A,B)$ and $[[A,B]]$ coincide \cite{Kasparov}, therefore there is a
one-to-one correspondence between homotopy classes of genuine and
asymptotic homomorphisms into the Calkin algebras $Q(B)$ and one has
$\Ext^{ph}(A,B)=0$.

}
\end{rmk}

\begin{prob}
{\rm
Does there exist a separable $C^*$-algebra $A$ such that
the $\Ext^{ph}(A,B)$ is non-zero for some $B$ ?

}
\end{prob}

Our definition of homotopy in $\Ext^{as}(A,B)$ is weaker than the homotopy
of asymptotic homomorphisms in $[[A,Q(B)]]$, so there is a surjective map
 \be\label{p}
p:[[A,Q(B)]]\ar\Ext^{as}(A,B).
 \ee
It would be interesting to compare the composition $\til{CH}\circ p$ with
the would-be boundary map $[[A,Q(B)]]\ar[[SA,B]]$ which would exist if the
exact sequences of the $E$-theory could be generalized to the
non-separable short exact sequence $B\ar M(B)\ar Q(B)$.

\begin{prob}
{\rm
Is the map $p$ (\ref{p}) injective ?

}
\end{prob}

%%%%%%%%%%%%%%%%%%%%%%%%%%%%%%%%%%%%%%%%%%%%%%%%%%%%%%%%%%%%%%%%%%%

\vspace{2cm}
\noindent
V.~M.~Manuilov\\
Dept. of Mech. and Math.,\\
Moscow State University,\\
Moscow, 119899, RUSSIA\\
e-mail: manuilov@mech.math.msu.su

\end{document}